\documentclass[12pt]{article}
\usepackage{amsfonts,amsmath,amsxtra}
\usepackage{latexsym}
\usepackage{amssymb}

\usepackage{comment}

\def\hybrid{\topmargin 0pt      \oddsidemargin 0pt
        \headheight 0pt \headsep 0pt
        \textwidth 16.5cm
        \textheight 23cm
        \marginparwidth 0.0in
        \parskip 5pt plus 1pt   \jot = 1.5ex}
\catcode`\@=11
\def\marginnote#1{}
\newcount\hour
\newcount\minute
\newtoks\amorpm
\hour=\time\divide\hour by60 \minute=\time{\multiply\hour by60
\global\advance\minute by-\hour}
\edef\standardtime{{\ifnum\hour<12 \global\amorpm={am}%
        \else\global\amorpm={pm}\advance\hour by-12 \fi
        \ifnum\hour=0 \hour=12 \fi
      \number\hour:\ifnum\minute<10 0\fi\number\minute\the\amorpm}}
\edef\militarytime{\number\hour:\ifnum\minute<10 0\fi\number\minute}

\def\draftlabel#1{{\@bsphack\if@filesw {\let\thepage\relax
   \xdef\@gtempa{\write\@auxout{\string
      \newlabel{#1}{{\@currentlabel}{\thepage}}}}}\@gtempa
   \if@nobreak \ifvmode\nobreak\fi\fi\fi\@esphack}
        \gdef\@eqnlabel{#1}}
\def\@eqnlabel{}
\def\@vacuum{}
\def\draftmarginnote#1{\marginpar{\raggedright\scriptsize\tt#1}}

\def\draft{\oddsidemargin -0.1truein
        \def\@oddfoot{\sl preliminary draft \hfil
        \rm\thepage\hfil\sl\today\quad\militarytime}
        \let\@evenfoot\@oddfoot \overfullrule 3pt
        \let\label=\draftlabel
        \let\marginnote=\draftmarginnote
\def\@eqnnum{{\rm (\theequation)}
\rlap{\kern\marginparsep\tt\@eqnlabel}%
\global\let\@eqnlabel\@vacuum}  }

\newfont{\Bbbb}{msbm7 scaled 1\@ptsize00}
\newcommand{\zs}{\raise-1pt\hbox{$\mbox{\Bbbb Z}$}}

scaled 1\@ptsize00

\font\sevenmsa=msam6 
\newfam\msafam
\textfont\msafam=\sevenmsa
\def\hexnumber@#1{\ifnum#1<10 \number#1\else
\ifnum#1=10 A\else\ifnum#1=11 B\else\ifnum#1=12 C\else \ifnum#1=13
D\else\ifnum#1=14 E\else\ifnum#1=15 F\fi\fi\fi\fi\fi\fi\fi}
\def\msa@{\hexnumber@\msafam}
\def\llcorner{\delimiter"4\msa@78\msa@78 }
\def\lrcorner{\delimiter"5\msa@79\msa@79 }
\mathchardef\blacktriangleright="3\msa@49
\mathchardef\blacktriangleleft="3\msa@4A \font\tenmsb=msbm10 scaled
1\@ptsize00
\newfam\msbfam
\textfont\msbfam=\tenmsb \scriptfont\msbfam=\tenmsb

\newdimen\Squaresize \Squaresize=14pt
\newdimen\Thickness \Thickness=0.5pt

\def\Square#1{\hbox{\vrule width \Thickness
   \vbox to \Squaresize{\hrule height \Thickness\vss
      \hbox to \Squaresize{\hss#1\hss}
   \vss\hrule height\Thickness}
\unskip\vrule width \Thickness} \kern-\Thickness}

\def\Vsquare#1{\vbox{\Square{$#1$}}\kern-\Thickness}

\def\numberbysection{\@addtoreset{equation}{section}
        \def\theequation{\thesection.\arabic{equation}}}
\numberbysection

\renewcommand{\theequation}{\thesection.\arabic{equation}}
\def\titlepage{\@restonecolfalse\if@twocolumn\@restonecoltrue\onecolumn
     \else \newpage \fi \thispagestyle{empty}\c@page\z@
        \def\thefootnote{\fnsymbol{footnote}} }

\def\endtitlepage{\if@restonecol\twocolumn \else  \fi
        \def\thefootnote{\arabic{footnote}}
        \setcounter{footnote}{0}}  
\relax

\hybrid
\makeatletter
\newdimen\normalarrayskip            
\newdimen\minarrayskip               
\normalarrayskip\baselineskip \minarrayskip\jot
\newif\ifold             \oldtrue            \def\new{\oldfalse}
\def\arraymode{\ifold\relax\else\displaystyle\fi}
\def\eqnumphantom{\phantom{(\theequation)}} 
\def\@arrayskip{\ifold\baselineskip\z@\lineskip\z@
     \else
     \baselineskip\minarrayskip\lineskip1\baselineskip\fi}

\def\@arrayclassz{\ifcase \@lastchclass \@acolampacol \or
\@ampacol \or \or \or \@addamp \or
   \@acolampacol \or \@firstampfalse \@acol \fi
\edef\@preamble{\@preamble
  \ifcase \@chnum
     \hfil$\relax\arraymode\@sharp$\hfil
     \or $\relax\arraymode\@sharp$\hfil
     \or \hfil$\relax\arraymode\@sharp$\fi}}

\def\@array[#1]#2{\setbox\@arstrutbox=\hbox{\vrule
     height\arraystretch \ht\strutbox
     depth\arraystretch \dp\strutbox
width\z@}\@mkpream{#2}\edef\@preamble{\halign \noexpand\@halignto
\bgroup \tabskip\z@ \@arstrut \@preamble \tabskip\z@ \cr}%
\let\@startpbox\@@startpbox \let\@endpbox\@@endpbox
  \if #1t\vtop \else \if#1b\vbox \else \vcenter \fi\fi
  \bgroup \let\par\relax
  \let\@sharp##\let\protect\relax
  \@arrayskip\@preamble}
\def\eqnarray{\stepcounter{equation}%
              \let\@currentlabel=\theequation
              \global\@eqnswtrue
              \global\@eqcnt\z@
              \tabskip\@centering              
              \let\\=\@eqncr
              $$%
            \halign to \displaywidth  \bgroup
             \eqnumphantom \@eqnsel
      \hskip\@centering                               
    $\displaystyle  \tabskip\z@ {##}$%
    &\global\@eqcnt\@ne \hskip 2\arraycolsep
         $ \displaystyle  \arraymode{##}$\hfil
    &\global\@eqcnt\tw@ \hskip 2\arraycolsep
         $\displaystyle\tabskip\z@{##}$\hfil
         \tabskip\@centering
    &{##}\tabskip\z@\cr}
\makeatother

\newtheorem{te}{Theorem}[section]
\newtheorem{de}{Definition}[section]
\newtheorem{cor}{Corollary}[section]

\newtheorem{rem}{Remark}[section]

\newcommand\bqa{\begin{eqnarray}}
\newcommand\eqa{\end{eqnarray}}
\def\be{\begin{eqnarray}\new\begin{array}{cc}}
\def\ee{\end{array}\end{eqnarray}}

\def\beq{\begin{equation}}
\def\eeq{\end{equation}}
\def\bse{\begin{subequations}}                
\def\ese{\end{subequations}}
\def\bp{\begin{pmatrix}}
\def\ep{\end{pmatrix}}

\def\proof{\noindent {\it Proof}. }

\def\stack#1#2{\raise0.7pt\hbox{$\mathrel{\mathop{#2}\limits^{#1}}$}}
\def\tr{\triangleright}
\def\tl{\triangleleft}
\def\sem{\mathsurround=0pt \raise1pt
\hbox{$\scriptscriptstyle>\!\!$}\:\!\!\tl}
\def\mes{\mathsurround=0pt \tr\!\:\!\raise0.8pt
\hbox{$\scriptscriptstyle\!\!<$}\,}
\def\]{\mathsurround=0pt ]\raise-2pt\hbox{$_\ast$}}

\def\<{\langle}
\def\>{\rangle}

\def\we{\raise-1pt\hbox{$\,\stackrel{\wedge}{,}\,$}}
\def\tr{{\rm tr}\,}

\newcounter{pac}[section]

\newcounter{pacc}[subsection]

0
\begin{document}

\setcounter{pac}{0}
\setcounter{footnote}0

\begin{center}

\phantom.
\bigskip
{\Large\bf On modular double of semisimple quantum groups}

\vspace{1cm}

\bigskip\bigskip
{\large  Pavel Sultanich \footnote {E-mail:  sultanichp@gmail.com}},\\
\bigskip
{\it Moscow Center for Continuous Mathematical Education, 119002, Bolshoy Vlasyevsky Pereulok 11, Moscow, Russia
}\\
\bigskip

\end{center}




\begin{abstract}
\noindent

In this note we propose a construction of the Hopf  algebra of a complex analog of devided powers of  the Weyl generators  of a semisimple simply-laced quantum group. Here  we consider  the generators as
positive, self-adjoint operators. In particular, we generalize the Lusztig relations \cite{Lu 1} on the usual divided powers of generators of a quantum group  to the case of complex devided powers of  generators. These relations, some of which were known \cite{Ip} present the complete set of defining relations.  As a by-product result, the pure algebraic definition of the Faddeev modular double  in the case of semisimple simply-laced quantum groups is formulated. Finally, we introduce an infinite dimension  version of the Gelfand-Zetlin finite-dimensional
representation of the modular double $M_{q,\tilde{q}}(\mathfrak{gl}(N))$.

\end{abstract}
\vspace{5 mm}

\section{Introduction}

The notion of modular double of a quantum group was originally introduced by Faddeev \cite{F2}. In case  $U_{q}(\mathfrak{sl}(2,\mathbb{R}))$, $q = e^{\pi\imath b^{2}}$, $b^{2}\in \mathbb{R}\setminus\mathbb{Q}$ he noticed that the class of representations of modular double, analogous to the principal series representations of $SL_2(R)$ has a remarkable duality, similar to the duality of non-commutative tori, discovered in \cite{Rif}. More generally, modular double plays an important role in Liouville theory \cite{PT}, \cite{FKV}, relativistic Toda model \cite{KLSTS} and some other problems of mathematical physics. Recently a considerable progress has been made in generalization of modular double to the case of $U_{q}(\mathfrak{g})$, where $\mathfrak{g}$ is a semisimple Lie algebra of rank $r$ \cite{FrIp}, \cite{Ip}, \cite{Ip2}. In these papers certain powers of generators have been used in order to construct the dual quantum group. This idea was first introduced in the case of $U_{q}(\mathfrak{sl}(2,\mathbb{R}))$ in \cite{BT}.

Let $U_{q}(\mathfrak{g})$ be the quantum group corresponding to semisimple simply-laced Lie algebra $\mathfrak{g}$ of rank $r$ with Cartan matrix $a_{ij} = a_{ji} \in\{0, -1\}$, for $i\ne j$, $1\le i,j\le r$. Let us assume that generators of $U_{q}(\mathfrak{g})$ are positive self-adjoint operators, then their complex powers can be defined \cite{Sh book}. In Theorem 4.1 the Hopf algebra of arbitrary complex
analogs of devided powers of generators has been constructed, which is a generalization of the algebra of divided powers $\frac{X^{N}}{[N]_{q}!}$ studied by Lusztig \cite{Lu 1}.

The Hopf algebra of complex  devided powers is generated by elements $K_{j}^{\imath p}$, $\mathcal{E}_{j}^{(\imath s)}$, $\mathcal{F}_{j}^{(\imath t)}$, $1\le j\le r$ subjected to certain integral relations. After specification $\imath s = N$, $\imath t = M$, $\imath p = L$, where $N$, $M$, $L$ are positive integers, all the defining integral relations are reduced to the identities for the devided powers of generators of quantum group, proved in \cite{Lu 1}. Some of these relations were known \cite{Ip}, however the other part of the relations including the generalized Kac identities is novel (see Theorem 3.1).

The Hopf algebra of complex devided powers is constructed using the classical Drinfeld's double technique \cite{D}. This Hopf algebra contains an important subalgebra generated by $K_{j}$, $\mathcal{E}_{j}$, $\mathcal{F}_{j}$ and $\tilde{K}_{j} = K_{j}^{b^{-2}}$, $\tilde{\mathcal{E}}_{j} = \mathcal{E}_{j}^{b^{-2}}$, $\tilde{\mathcal{F}}_{j} = \mathcal{F}_{j}^{b^{-2}}$, $1\le j\le r$. The integral relations of the full Hopf algebra of complex devided powers are reduced to the standard relations of quantum groups $U_{q}(\mathfrak{g})$, generated by $K_{j}$, $\mathcal{E}_{j}$, $\mathcal{F}_{j}$, $1\le j\le r$ and $U_{\tilde{q}}(\mathfrak{g})$ generated by $\tilde{K}_{j}$, $\tilde{\mathcal{E}}_{j}$, $\tilde{\mathcal{F}}_{j}$, $1\le j\le r$ with $\tilde{q} = e^{\pi\imath b^{-2}}$, together with non-trivial cross-relations. This subalgebra has been the subject of active research recently \cite{BT}, \cite{FrIp}, \cite{Ip}, \cite{Ip2}, and it is usually associated with the modular double construction. However, the phenomenon of non-triviality of cross-relations was studied only for special representations. In the present paper we propose a natural derivation of full set of relations of modular double from the integral relations of the algebra of complex powers (see Theorem 5.1).

Next we study an algebraic definition of modular double. The modular double is a Hopf algebra with the two sets of generators: generators $K_{j}$, $\mathcal{E}_{j}$, $\mathcal{F}_{j}$, $1\le j\le r$ subjected to standard relations of $U_{q}(\mathfrak{g})$, and additional generators $\tilde{K}_{j}$, $\tilde{\mathcal{E}}_{j}$, $\tilde{\mathcal{F}}_{j}$, $1\le j\le r$ subjected the standard relations of $U_{\tilde{q}}(\mathfrak{g})$,
 together with the  same cross relations but free of the
  constraints  $\tilde{K}_{j} = K_{j}^{b^{-2}}$, $\tilde{\mathcal{E}}_{j} = \mathcal{E}_{j}^{b^{-2}}$, $\tilde{\mathcal{F}}_{j} = \mathcal{F}_{j}^{b^{-2}}$, $1\le j\le r$. The algebraic definition is motivated by the fact that these Hopf algebras have two types of realizations of the principal series representations. In Theorem 5.2 a realization of representation of the modular double of $U_{q}(\mathfrak{gl}(N))$ of the first type is given.  This realization of the principal series representation \cite{GKL1}, \cite{GKL2}, \cite{GKLO} of the modular double can be viewed as an infinite-dimensional analog of the Gelfand-Zetlin finite-dimensional representation for classical groups \cite{GZ}.


The realizations of second type \cite{FrIp}, \cite{Ip2} are the principal series representations in Lusztig parametrization. These representations are $q$-deformed versions of representations of universal enveloping algebra $U(\mathfrak{g})$ by differential operators in Lusztig parametrization introduced in \cite{GLO}(sections 2.4.1-2.4.4) for classical series Lie algebras. The non-trivial properties of the representations are transcendental relations $\tilde{K}_{j} = K_{j}^{b^{-2}}$, $\tilde{\mathcal{E}}_{j} = \mathcal{E}_{j}^{b^{-2}}$, $\tilde{\mathcal{F}}_{j} = \mathcal{F}_{j}^{b^{-2}}$, $1\le j\le r$. Using the classical limit of this realization leads to construction of Whittaker function in terms of the stationary phase integral, generalizing Givental's formula, see \cite{GLO}.

{\bf Acknowledgements:} The research was supported by    RSF  (project № 16-11-10075). I am grateful to D.R. Lebedev for the statement  of the problem and his  interest  in this work.

\section{Preliminaries}

We start with the definition of quantum groups following \cite{ChPr},\cite{Lu book}.
Let $(a_{ij})_{1\le i,j\le r}$ be Cartan matrix of semisimple Lie algebra $\mathfrak{g}$ of rank $r$. Let $\mathfrak{b}_{\pm}\subset \mathfrak{g}$ be opposite Borel subalgebras. For simplicity let us restrict ourselves to the simply-laced case $a_{ii} = 2$, $a_{ij} = a_{ji} = \{0,-1\}$, $i\ne j$. Let $U_{q}(\mathfrak{g})$ $(q = e^{\pi\imath b^{2}}$, $b^{2}\in \mathbb{R}\setminus \mathbb{Q})$ be the quantum group with generators $E_{j}$, $F_{j}$, $K_{j} = q^{H_{j}}$, $1\le j \le r$ and relations
\begin{equation}
    K_{i}K_{j} = K_{j}K_{i},
\end{equation}
\begin{equation}
    K_{i}E_{j} = q^{a_{ij}}E_{j}K_{i},
\end{equation}
\begin{equation}
    K_{i}F_{j} = q^{-a_{ij}}F_{j}K_{i},
\end{equation}
\begin{equation}
    E_{i}F_{j} - F_{j}E_{i} = \delta_{ij}\frac{K_{i} - K_{i}^{-1}}{q-q^{-1}}.
\end{equation}
For $a_{ij} = 0$ we have
\begin{equation}
    E_{i}E_{j} = E_{j}E_{i},
\end{equation}
\begin{equation}
    F_{i}F_{j} = F_{j}F_{i}.
\end{equation}
For $a_{ij} = -1$ we have
\begin{equation}
    E_{i}^{2}E_{j} - (q+q^{-1})E_{i}E_{j}E_{i} + E_{j}E_{i}^{2} = 0,
\end{equation}
\begin{equation}
    F_{i}^{2}F_{j} - (q+q^{-1})F_{i}F_{j}F_{i} + F_{j}F_{i}^{2} = 0,
\end{equation}
Coproduct is given by
\begin{equation}
\Delta E_{j} = E_{j}\otimes 1 + K_{j}^{-1}\otimes E_{j},
\end{equation}
\begin{equation}
\Delta F_{j} = 1\otimes F_{j} + F_{j}\otimes K_{j},
\end{equation}
\begin{equation}
\Delta K_{j} = K_{j}\otimes K_{j}.
\end{equation}

Non-compact quantum dilogarithm $G_{b}(z)$ is a special function introduced in \cite{F1} (see also \cite{F0}, \cite{FKV}, \cite{V}, \cite{Ka}, \cite{KLSTS}, \cite{BT}). It is defined as follows
\begin{equation}
\log G_{b}(z) = \log\bar{\zeta}_{b} - \int\limits_{\mathbb{R}+\imath 0} \frac{dt}{t}\frac{e^{zt}}{(1-e^{bt})(1-e^{b^{-1}t})},
\end{equation}
where $Q = b+b^{-1}$ and $\zeta_{b} = e^{\frac{\pi\imath}{4} + \frac{\pi\imath(b^{2}+b^{-2})}{12}}$. Note, that $G_{b}(z)$ is closely related to the double sine function $S_{2}(z|\omega_{1},\omega_{2})$, see eq.(A.22) in \cite{KLSTS}.

Below we outline some properties of $G_{b}(z)$, for details see appendix.\\*
1. The function $G_{b}(z)$ has simple poles and zeros at the points
\begin{equation}
    z = -n_{1}b -n_{2}b^{-1},
\end{equation}
\begin{equation}
    z = Q +n_{1}b + n_{2}b^{-1},
\end{equation}
respectively, where $n_{1}$,$n_{2}$ are nonnegative integer numbers.\\*
2. $G_{b}(z)$ has the following asymptotic behavior:
\begin{equation}
 G_{b}(z) \sim
 \begin{cases} \bar{\zeta}_{b}, Im z \rightarrow +\infty ,\\ \zeta_{b} e^{\pi\imath z(z-Q)}, Im z \rightarrow -\infty . \end{cases}
\end{equation}
3. Functional equation:
\begin{equation}
G_{b}(z +b^{\pm 1}) = (1-e^{2\pi\imath b^{\pm 1}z})G_{b}(z).
\end{equation}
4. Reflection formula:
\begin{equation}
G_{b}(z)G_{b}(Q-z) = e^{\pi\imath z(z-Q)}.
\end{equation}

For $1\le j \le r$ let us introduce the following rescaled generators
\begin{equation}\label{rescaled E}
    \mathcal{E}_{j} = -\imath (q-q^{-1})E_{j},
\end{equation}
\begin{equation}\label{rescaled F}
    \mathcal{F}_{j} = -\imath (q-q^{-1})F_{j}.
\end{equation}
We will assume the elements $\mathcal{E}_{j}$, $\mathcal{F}_{j}$, $K_{j}$, $1 \le i \le r$ to be represented by positive self-adjoint operators. For such operator $A$ one can define its arbitrary complex powers \cite{Sh book} so that the following properties are satisfied \\*
1. $A^{\imath s}$ is holomorphic in the half-plane $Im(s)> -k$ for any $k\in \mathbb{Z}$. \\*
2. Group property
\begin{equation}
    A^{\imath s_{1}}A^{\imath s_{2}} = A^{\imath s_{1} + \imath s_{2}}.
\end{equation}
3. For $\imath s = k$, $k\in\mathbb{Z}$, $A^{k}$ is an ordinary power of $A$.\\*
Next, we define the arbitrary devided powers of $A$ by
\begin{equation}\label{complex devided power}
    A^{(\imath s)} = G_{b}(-\imath bs)A^{\imath s}.
\end{equation}
We are going to consider the algebra spanned by the elements $\mathcal{E}_{j}^{(\imath s)}$, $\mathcal{F}_{j}^{(\imath t)}$, $K_{j}^{\imath p}$, $1 \le j \le r$.

\section{Generalized Kac's identity}

Below we give a generalization of the Kac's identity \cite{Lu 1}, eq.(4.1a) to the case of complex powers of generators.

\begin{te}
Let $q = e^{\pi\imath b^{2}}$, $(b^{2}\in \mathbb{R}\setminus \mathbb{Q})$ and  let $K_{j} = q^{H_{j}}$, $\mathcal{E}_{j} = -\imath (q-q^{-1})E_{j}$, $\mathcal{F}_{j} = -\imath (q-q^{-1})F_{j}$, $1\le j\le r$ be positive self-adjoint operators and let $A^{(\imath s)}$ be defined by (\ref{complex devided power}). Then the following generalized Kac's identity holds:
\bigskip
\begin{equation}\begin{split}
\mathcal{E}_{j}^{(\imath s)}\mathcal{F}_{j}^{(\imath t)} = \int\limits_{\mathcal{C}} d\tau e^{\pi bQ\tau}\mathcal{F}_{j}^{(\imath t+\imath\tau)}K_{j}^{-\imath \tau}
\frac{G_{b}(\imath b\tau)G_{b}(-bH_{j} + \imath b(s+t+\tau))}{G_{b}(-bH_{j}+\imath b(s+t+2\tau))}\mathcal{E}_{j}^{(\imath s + \imath \tau)},
\end{split}\end{equation}
\bigskip
where the contour $\mathcal{C}$ goes along the real axis above the sequences of poles going down: \\* $\tau = -s-\imath n_{1}-\imath n_{2}b^{-2}$, $\tau = -t-\imath n_{1}-\imath n_{2}b^{-2}$, $\tau = -\frac{\imath b^{-1}Q}{2}-\frac{\imath H_{j}}{2} -\frac{s}{2}-\frac{t}{2}-\frac{\imath n_{1}}{2}-\frac{\imath b^{-2}n_{2}}{2}$, \\* and below the sequences of poles going up:\\*
$\tau = \imath n_{1} +\imath n_{2}b^{-2}$, $\tau = -\imath H_{j}-s-t+\imath n_{1}+\imath n_{2}b^{-2}$,\\*
where $n_{1}$, $n_{2}$ are non-negative integers.
\end{te}
$\proof$
The proof is based on the Drinfeld's double construction \cite{D} and will be published elsewhere.
$\Box$

Let $N$, $M$ be positive integers. Define the integer devided powers of generators $E_{j}$, $F_{j}$, $1\le j \le r$ of quantum group $U_{q}(\mathfrak{g})$:
\begin{equation}
E_{j}^{(N)} = \frac{\prod\limits_{k=1}^{N}(q-q^{-1})}{\prod\limits_{k=1}^{N}(q^{k}-q^{-k})}E_{j}^{N},
\end{equation}
\begin{equation}
F_{j}^{(M)} = \frac{\prod\limits_{k=1}^{M}(q-q^{-1})}{\prod\limits_{k=1}^{M}(q^{k}-q^{-k})}F_{j}^{M}.
\end{equation}
\begin{cor}
Let $\imath s = N$, $\imath t = M$, where $N$, $M$ are non-negative integers. Then the generalized Kac's identity reduces to the standard Kac's identity \cite{Lu 1}, eq.(4.1a):
\begin{equation}
    E_{j}^{(N)}F_{j}^{(M)} =
\sum\limits_{n=0}^{min(N,M)} F_{j}^{(M-n)}
\frac{\prod\limits_{k=1}^{n}(q^{-N-M+n+k}K_{j}-q^{N+M-n-k}K_{j}^{-1})}{\prod\limits_{k=1}^{n}(q^{k}-q^{-k})}
E_{j}^{(N-n)},
\end{equation}
where $1\le j\le r$. Substituting into this formula $N = M = 1$ one recovers the relation
\begin{equation}
    [E_{j},F_{j}] = \frac{K_{j}-K_{j}^{-1}}{q-q^{-1}},
\end{equation}
for $1\le j\le r$.
\end{cor}
$\proof$
Rewrite the generalized Kac's identity in the following form:
$$
\mathcal{E}_{j}^{\imath s}\mathcal{F}_{j}^{\imath t} = \int d\tau e^{\pi bQ\tau} \frac{G_{b}(\imath b\tau)G_{b}(-\imath bt-\imath b\tau)}{G_{b}(-\imath bt)}
\frac{G_{b}(-\imath bs-\imath b\tau)}{G_{b}(-\imath bs)}\mathcal{F}_{j}^{\imath t+\imath\tau} K_{j}^{-\imath\tau}
\frac{G_{b}(-bH_{j}+\imath b(s+t+\tau))}{G_{b}(-bH_{j}+\imath b(s+t+2\tau))}\mathcal{E}_{j}^{\imath s+\imath\tau}.
$$
Set $\imath t  = M$ to be some nonnegative integer. Then using delta-distribution formula \cite{Ip}
\begin{equation}\begin{split}
    \frac{G_{b}(x)G_{b}(-N_{1}b-N_{2}b^{-1}-x)}{G_{b}(-N_{1}b-N_{2}b^{-1})} =
    \sum\limits_{n_{1}=0}^{N_{1}}\sum\limits_{n_{2}=0}^{N_{2}}\frac{\prod\limits_{k_{1}=1}^{N_{1}}(1-q^{-2k_{1}})}{\prod\limits_{k_{1}=1}^{n_{1}}(1-q^{-2k_{1}})\prod\limits_{k_{1}=1}^{N_{1}-n_{1}}(1-q^{-2k_{1}})} \\
   \times \frac{\prod\limits_{k_{2}=1}^{N_{2}}(1-\tilde{q}^{-2k_{2}})}{\prod\limits_{k_{2}=1}^{n_{2}}(1-\tilde{q}^{-2k_{2}})\prod\limits_{k_{2}=1}^{N_{2}-n_{2}}(1-\tilde{q}^{-2k_{2}})}
    \delta(x+n_{1}b+n_{2}b^{-1}),
\end{split}\end{equation}
we obtain
$$
\mathcal{E}_{j}^{\imath s}\mathcal{F}_{j}^{M}  = \sum\limits_{n=0}^{M}
\frac{\prod\limits_{k=1}^{M}(1-q^{-2k})}{\prod\limits_{k=1}^{n}(1-q^{-2k})\prod\limits_{k=1}^{M-n}(1-q^{-2k})}\times
$$
$$
\int d\tau e^{\pi bQ\tau}\delta(\imath b\tau + nb)
\frac{G_{b}(-\imath bs-\imath b\tau)}{G_{b}(-\imath bs)}\mathcal{F}_{j}^{M+\imath\tau} K_{j}^{-\imath\tau}
\frac{G_{b}(-bH_{j}+\imath bs+Mb+\imath b\tau)}{G_{b}(-bH_{j}+\imath bs+Mb+2\imath b\tau)}\mathcal{E}_{j}^{\imath s+\imath\tau} =
$$

$$
\sum\limits_{n=0}^{M}
\frac{\prod\limits_{k=1}^{M}(1-q^{-2k})}{\prod\limits_{k=1}^{n}(1-q^{-2k})\prod\limits_{k=1}^{M-n}(1-q^{-2k})}
(-1)^{n}q^{n}\times
$$
$$
\frac{G_{b}(-\imath bs+nb)}{G_{b}(-\imath bs)}\mathcal{F}_{j}^{M-n}K^{n}
\frac{G_{b}(-bH_{j}+\imath bs+Mb-2nb+nb)}{G_{b}(-bH_{j}+\imath bs+Mb-2nb)}\mathcal{E}_{j}^{\imath s-n}  =
$$
$$
\sum\limits_{n=0}^{M}
\frac{\prod\limits_{k=1}^{M}(1-q^{-2k})}{\prod\limits_{k=1}^{n}(1-q^{-2k})\prod\limits_{k=1}^{M-n}(1-q^{-2k})}
(-1)^{n}q^{n}
\prod\limits_{k=0}^{n-1}(1-q^{2k}e^{2\pi b^{2}s})\mathcal{F}_{j}^{M-n}K^{n}\times
$$
$$
\times\prod\limits_{k=0}^{n-1}(1-q^{2k}e^{2\pi\imath b(-bH_{j}+\imath bs+Mb-2nb)})\mathcal{E}_{j}^{\imath s-n}.
$$
To write the last line we have used the functional equation for quantum dilogarithm:
\begin{equation}
\frac{G_{b}(x+n_{1}b+n_{2}b^{-1})}{G_{b}(x)} =
\prod\limits_{k_{1} =0}^{n_{1}-1}(1-q^{2k_{1}}e^{2\pi\imath bx})\prod\limits_{k_{2} =0}^{n_{2}-1}(1-\tilde{q}^{2k_{2}}e^{2\pi\imath b^{-1}x}).
\end{equation}
We have
$$
\mathcal{E}_{j}^{\imath s}\mathcal{F}_{j}^{M} =
\sum\limits_{n=0}^{M}
\frac{\prod\limits_{k=1}^{M}(q^{k}-q^{-k})}{\prod\limits_{k=1}^{n}(q^{k}-q^{-k})\prod\limits_{k=1}^{M-n}(q^{k}-q^{-k})}(-1)^{n}
\prod\limits_{k=0}^{n-1}(q^{-k}e^{-\pi b^{2}s} - q^{k}e^{\pi b^{2}s})\times
$$
$$
\mathcal{F}_{j}^{M-n}
\prod\limits_{k=0}^{n-1}(e^{\pi b^{2}s-\pi\imath b^{2}M+2\pi\imath b^{2}n-\pi\imath b^{2}k}K_{j} -e^{-\pi b^{2}s+\pi\imath b^{2}M-2\pi\imath b^{2}n+\pi\imath b^{2}k}K_{j}^{-1})\mathcal{E}_{j}^{\imath s-n}.
$$
If $M=1$ then we get the fomula from \cite{BT}, eq.(3.25)
\begin{equation}
[\mathcal{E}_{j}^{\imath s},\mathcal{F}_{j}] = -(q^{\imath s} - q^{-\imath s})(q^{H_{j}-\imath s+1} - q^{-H_{j}+\imath s-1})\mathcal{E}_{j}^{\imath s-1},
\end{equation}
which was derived in representation in \cite{BT}.

Now let $M$ again be arbitrary positive integer, set $\imath s = N$ and substitute
$$
\mathcal{E}_{j}^{N} = (-\imath)^{N}\prod\limits_{k=1}^{N}(q^{k}-q^{-k})E_{j}^{(N)},
$$
$$
\mathcal{F}_{j}^{M} = (-\imath)^{M}\prod\limits_{k=1}^{M}(q^{k}-q^{-k})F_{j}^{(M)},
$$
to obtain the Kac formula:
$$
E_{j}^{(N)}F_{j}^{(M)} =
\sum\limits_{n=0}^{min(N,M)} F_{j}^{(M-n)}
\frac{\prod\limits_{k=1}^{n}(q^{-N-M+n+k}K_{j}-q^{N+M-n-k}K_{j}^{-1})}{\prod\limits_{k=1}^{n}(q^{k}-q^{-k})}
E_{j}^{(N-n)}.
$$

Let in the generalized Kac's identity $\imath s = 1$ and let $\imath t$ be arbitrary. Then by evaluating the integral as above one can prove the formula \cite{BT}, eq.(3.25)
\begin{equation}
    [\mathcal{E}_{j},\mathcal{F}_{j}^{\imath t}] = -(q^{\imath t}-q^{-\imath t})\mathcal{F}_{j}^{\imath t-1}(q^{H_{j}-\imath t+1} - q^{-H_{j}+\imath t-1}).
\end{equation}

$\Box$

Let $\tilde{q} = e^{\pi\imath b^{-2}}$. Define the dual generators $\tilde{K}_{j}$, $\tilde{E}_{j}$, $\tilde{F}_{j}$, $1\le j\le r$ by
\begin{equation}
    \tilde{K}_{j} = K_{j}^{b^{-2}},
\end{equation}
\begin{equation}
    \tilde{E}_{j} = \frac{\imath}{\tilde{q}-\tilde{q}^{-1}}\mathcal{E}_{j}^{b^{-2}},
\end{equation}
\begin{equation}
    \tilde{F}_{j} = \frac{\imath}{\tilde{q}-\tilde{q}^{-1}}\mathcal{F}_{j}^{b^{-2}},
\end{equation}
and devided powers of $\tilde{E}_{j}$ and $\tilde{F}_{j}$, $1\le j\le r$:
\begin{equation}
\tilde{E}_{j}^{(N)} = \frac{\prod\limits_{k=1}^{N}(\tilde{q}-\tilde{q}^{-1})}{\prod\limits_{k=1}^{N}(\tilde{q}^{k}-\tilde{q}^{-k})}\tilde{E}_{j}^{N},
\end{equation}
\begin{equation}
\tilde{F}_{j}^{(M)} = \frac{\prod\limits_{k=1}^{M}(\tilde{q}-\tilde{q}^{-1})}{\prod\limits_{k=1}^{M}(\tilde{q}^{k}-\tilde{q}^{-k})}\tilde{F}_{j}^{M}.
\end{equation}

\begin{cor}
By setting $\imath s = Nb^{-2}$, $\imath t = Mb^{-2}$ in the generalized Kac's identity, where $N$, $M$ are non-negative integers, one obtains the standard Kac's identity for the dual quantum group:
\begin{equation}
    \tilde{E}_{j}^{(N)}\tilde{F}_{j}^{(M)} =
\sum\limits_{n=0}^{min(N,M)} \tilde{F}_{j}^{(M-n)}
\frac{\prod\limits_{k=1}^{n}(\tilde{q}^{-N-M+n+k}\tilde{K}_{j}-\tilde{q}^{N+M-n-k}\tilde{K}_{j}^{-1})}{\prod\limits_{k=1}^{n}(\tilde{q}^{k}-\tilde{q}^{-k})}
\tilde{E}_{j}^{(N-n)},
\end{equation}
where $1\le j\le r$. Taking $N = M = 1$ we obtain one of the defining relations of the dual quantum group
\begin{equation}
[\tilde{E}_{j},\tilde{F}_{j}] = \frac{\tilde{K}_{j}-\tilde{K}_{j}^{-1}}{\tilde{q}-\tilde{q}^{-1}},
\end{equation}
$1\le j\le r$.
\end{cor}
$\proof$
Following the same steps as in the proof of previous corollary with the only difference that now $\imath t = Mb^{-2}$ and $\imath s = Nb^{-2}$ and substituting for $1\le j\le r$
$$
\mathcal{E}_{j}^{Nb^{-2}} = (-\imath)^{N}\prod\limits_{k=1}^{N}(\tilde{q}^{k}-\tilde{q}^{-k})\tilde{E}_{j}^{(N)},
$$
$$
\mathcal{F}_{j}^{Mb^{-2}} = (-\imath)^{M}\prod\limits_{k=1}^{M}(\tilde{q}^{k}-\tilde{q}^{-k})\tilde{F}_{j}^{(M)},
$$
we obtain the Kac's identity for the modular dual group with the parameter $\tilde{q} = e^{\pi\imath b^{-2}}$
$$
\tilde{E}_{j}^{(N)}\tilde{F}_{j}^{(M)} =
\sum\limits_{n=0}^{min(N,M)} \tilde{F}_{j}^{(M-n)}
\frac{\prod\limits_{k=1}^{n}(\tilde{q}^{-N-M+n+k}\tilde{K}_{j}-\tilde{q}^{N+M-n-k}\tilde{K}_{j}^{-1})}{\prod\limits_{k=1}^{n}(\tilde{q}^{k}-\tilde{q}^{-k})}
\tilde{E}_{j}^{(N-n)}.
$$

$\Box$

\begin{cor}
\begin{equation}
[\mathcal{E}_{j}^{\imath s},\mathcal{F}_{j}] = -(q^{\imath s} - q^{-\imath s})(q^{H_{j}-\imath s+1} - q^{-H_{j}+\imath s-1})\mathcal{E}_{j}^{\imath s-1},
\end{equation}
\begin{equation}
    [\mathcal{E}_{j},\mathcal{F}_{j}^{\imath t}] = -(q^{\imath t}-q^{-\imath t})\mathcal{F}_{j}^{\imath t-1}(q^{H_{j}-\imath t+1} - q^{-H_{j}+\imath t-1}),
\end{equation}
where $1\le j\le r$. Let $\tilde{\mathcal{E}}_{j} = \mathcal{E}_{j}^{b^{-2}}$, $\tilde{\mathcal{F}}_{j} = \mathcal{F}_{j}^{b^{-2}}$, $1\le j\le r$. Then setting $\imath s = b^{-2}$ and $\imath t = b^{-2}$ one obtains
\begin{equation}
[\tilde{\mathcal{E}}_{j},\mathcal{F}_{j}] = [\mathcal{E}_{j},\tilde{\mathcal{F}}_{j}] = 0.
\end{equation}
\end{cor}
$\proof$
See the proof of the Corollary 3.1.
$\Box$

\section{Hopf algebra of arbitrary complex devided powers}
Let $U_{q}(\mathfrak{g})$ be quantum group corresponding to semisimple simply-laced Lie algebra $\mathfrak{g}$ with Cartan matrix $(a_{ij})$, $1\le i,j\le r$.
Recall, that we have introduced the rescaled generators $\mathcal{E}_{j}$, $\mathcal{F}_{j}$ of $U_{q}(\mathfrak{g})$ by the formulas
\begin{equation}
    \mathcal{E}_{j} = -\imath (q-q^{-1})E_{j},
\end{equation}
\begin{equation}
    \mathcal{F}_{j} = -\imath (q-q^{-1})F_{j}.
\end{equation}
We have also defined arbitrary devided powers
\begin{equation}
    A^{(\imath s)} = G_{b}(-\imath bs)A^{\imath s}.
\end{equation}
For $a_{ij} = -1$ define non-simple root generators by
\begin{equation}
    \mathcal{E}_{ij} = \frac{q^{\frac{1}{2}}\mathcal{E}_{j}\mathcal{E}_{i}-q^{-\frac{1}{2}}\mathcal{E}_{i}\mathcal{E}_{j}}{q-q^{-1}},
\end{equation}
\begin{equation}
    \mathcal{F}_{ij} = \frac{q^{\frac{1}{2}}\mathcal{F}_{j}\mathcal{F}_{i}-q^{-\frac{1}{2}}\mathcal{F}_{i}\mathcal{F}_{j}}{q-q^{-1}}.
\end{equation}
The next theorem is a generalization of subalgebra of devided powers $\frac{X^{n}}{[n]_{q}!}$, $X\in U_{q}(\mathfrak{g})$ and relations \cite{Lu 1}, eq.(4.1 a) - eq.(4.1 j) to the entire Hopf algebra of complex devided powers of generators.
\begin{te}
The elements $K_{j}^{\imath p}$, $\mathcal{E}_{j}^{(\imath s)}$, $\mathcal{F}_{j}^{(\imath t)}$, $1\le j \le r$  generate associative and coassociative Hopf algebra ${A}(\mathfrak{g})$ with the following commutation relations and coproduct
\begin{equation}
    K_{i}^{\imath p_{1}}K_{j}^{\imath p_{2}} = K_{j}^{\imath p_{2}}K_{i}^{\imath p_{1}},
\end{equation}
\begin{equation}
K_{j}^{\imath p_{1}}K_{j}^{\imath p_{2}} = K_{j}^{\imath p_{1}+\imath p_{2}},
\end{equation}
\begin{equation}
\mathcal{E}_{j}^{(\imath s_{1})}\mathcal{E}_{j}^{(\imath s_{2})} = \frac{G_{b}(-\imath bs_{1})G_{b}(-\imath bs_{2})}{G_{b}(-\imath bs_{1}-\imath bs_{2})}\mathcal{E}_{j}^{(\imath s_{1}+\imath s_{2})} = \mathcal{E}_{j}^{(\imath s_{2})}\mathcal{E}_{j}^{(\imath s_{1})},
\end{equation}
\begin{equation}
\mathcal{F}_{j}^{(\imath t_{1})}\mathcal{F}_{j}^{(\imath t_{2})} = \frac{G_{b}(-\imath bt_{1})G_{b}(-\imath bt_{2})}{G_{b}(-\imath bt_{1}-\imath bt_{2})}\mathcal{F}_{j}^{(\imath t_{1}+\imath t_{2})} = \mathcal{F}_{j}^{(\imath t_{2})}\mathcal{F}_{j}^{(\imath t_{1})},
\end{equation}
\begin{equation}
K_{i}^{\imath p}\mathcal{E}_{j}^{(\imath s)} = e^{-\pi\imath b^{2}a_{ij}ps}\mathcal{E}_{j}^{(\imath s)}K_{i}^{\imath p},
\end{equation}
\begin{equation}
K_{i}^{\imath p}\mathcal{F}_{j}^{(\imath t)} = e^{\pi\imath b^{2}a_{ij}pt}\mathcal{F}_{j}^{(\imath t)}K_{i}^{\imath p}.
\end{equation}

If $a_{ij} = 0$, then
\begin{equation}
    \mathcal{E}_{i}^{(\imath s)}\mathcal{E}_{j}^{(\imath t)} = \mathcal{E}_{j}^{(\imath t)}\mathcal{E}_{i}^{(\imath s)},
\end{equation}
\begin{equation}
    \mathcal{F}_{i}^{(\imath s)}\mathcal{F}_{j}^{(\imath t)} = \mathcal{F}_{j}^{(\imath t)}\mathcal{F}_{i}^{(\imath s)}.
\end{equation}

For $a_{ij} = -1$ we have
\begin{equation} \label{generalized Serre E}
    \mathcal{E}_{i}^{(\imath s)}\mathcal{E}_{j}^{(\imath t)} = e^{-\pi\imath b^{2}st}\int\limits_{\Gamma_{1}} d\tau e^{\frac{\pi\imath b^{2}\tau^{2}}{2}-\pi bQ\tau}
    \mathcal{E}_{j}^{(\imath t-\imath\tau)}\mathcal{E}_{ij}^{(\imath\tau)}\mathcal{E}_{i}^{(\imath s-\imath\tau)},
\end{equation}
\begin{equation}\label{generalized Serre F}
    \mathcal{F}_{i}^{(\imath s)}\mathcal{F}_{j}^{(\imath t)} = e^{-\pi\imath b^{2}st}\int\limits_{\Gamma_{1}} d\tau e^{\frac{\pi\imath b^{2}\tau^{2}}{2}-\pi bQ\tau}
    \mathcal{F}_{j}^{(\imath t-\imath\tau)}\mathcal{F}_{ij}^{(\imath\tau)}\mathcal{F}_{i}^{(\imath s-\imath\tau)},
\end{equation}
where the contour $\Gamma_{1}$ goes above the pole at $\tau = 0$ and below the poles at $\tau = s$, $\tau = t$.\\*
For $i\ne j$ we have
\begin{equation}
    \mathcal{E}_{i}^{(\imath s)}\mathcal{F}_{j}^{(\imath t)} = \mathcal{F}_{j}^{(\imath t)}\mathcal{E}_{i}^{(\imath s)}.
\end{equation}

\begin{equation}\label{generalized Kac's identity}
\mathcal{E}_{j}^{(\imath s)}\mathcal{F}_{j}^{(\imath t)} = \int\limits_{\mathcal{C}} d\tau e^{\pi bQ\tau}\mathcal{F}_{j}^{(\imath t+\imath\tau)}K_{j}^{-\imath \tau}
\frac{G_{b}(\imath b\tau)G_{b}(-bH_{j} + \imath b(s+t+\tau))}{G_{b}(-bH_{j}+\imath b(s+t+2\tau))}\mathcal{E}_{j}^{(\imath s + \imath \tau)},
\end{equation}
where the contour $\mathcal{C}$ is defined in Theorem 3.1.\\*
Coproduct consistent with commutation relations has the form
\begin{equation}
\Delta K_{j}^{\imath p} = K_{j}^{\imath p}\otimes K_{j}^{\imath p},
\end{equation}
\begin{equation}
\Delta\mathcal{E}_{j}^{(\imath s)} = \int\limits_{\Gamma_{2}} d\tau \mathcal{E}_{j}^{(\imath s -\imath\tau)}K_{j}^{-\imath\tau}\otimes \mathcal{E}_{j}^{(\imath\tau)},
\end{equation}
where $\Gamma_{2}$ goes above the pole at $\tau = 0$ and below the pole at $\tau = s$.
\begin{equation}
\Delta\mathcal{F}_{j}^{(\imath t)} = \int\limits_{\Gamma_{3}} d\tau \mathcal{F}_{j}^{(\imath\tau)}\otimes \mathcal{F}_{j}^{(\imath t-\imath\tau)}K_{j}^{\imath\tau}.
\end{equation}
$\Gamma_{3}$ goes above the pole at $\tau = 0$ and below the pole at $\tau = t$.
\end{te}

\begin{rem}
The relations (\ref{generalized Serre E}), (\ref{generalized Serre F}) appeared in \cite{Ip}, eq.(6.16).
The formula (\ref{generalized Kac's identity}) is a new one.
\end{rem}

\begin{cor}
Let $\tilde{K}_{i} = K_{i}^{b^{-2}}$, $1\le i\le r$. Then any element from both sets $K_{i}$, $1\le i\le r$ and $\tilde{K}_{j}$, $1\le j\le r$ commute with any other element from these sets.
\end{cor}
$\proof$
The statement of the corollary follows straightforwardly from formulas (4.6), (4.7) if we take $\imath p_{1} = 1$, $\imath p_{2} = b^{-2}$.
$\Box$

\begin{cor}
Let $\tilde{\mathcal{E}}_{j} = \mathcal{E}_{j}^{b^{-2}}$ and $\tilde{\mathcal{F}}_{j} = \mathcal{F}_{j}^{b^{-2}}$, $1\le j\le r$. Then
\begin{equation}
    \tilde{\mathcal{E}}_{j}\mathcal{E}_{j} = \mathcal{E}_{j}\tilde{\mathcal{E}}_{j},
\end{equation}
\begin{equation}
    \tilde{\mathcal{F}}_{j}\mathcal{F}_{j} = \mathcal{F}_{j}\tilde{\mathcal{F}}_{j},
\end{equation}
where $1\le j\le r$.
\end{cor}
$\proof$
The statement of the corollary follows immediately from formulas (4.8), (4.9) in the case of $\imath s_{1} = 1$, $\imath s_{2} = b^{-2}$ and $\imath t_{1} = 1$, $\imath t_{2} = b^{-2}$.
$\Box$

\begin{cor}
The following relations hold
\begin{equation}
    K_{i}\mathcal{E}_{j} = q^{a_{ij}}\mathcal{E}_{j}K_{i},
\end{equation}
\begin{equation}
    K_{i}\mathcal{F}_{j} = q^{-a_{ij}}\mathcal{F}_{j}K_{i},
\end{equation}
\begin{equation}
    \tilde{K}_{i}\tilde{\mathcal{E}}_{j} = \tilde{q}^{a_{ij}}\tilde{\mathcal{E}}_{j}\tilde{K}_{i},
\end{equation}
\begin{equation}
    \tilde{K}_{i}\tilde{\mathcal{F}}_{j} = \tilde{q}^{-a_{ij}}\tilde{\mathcal{F}}_{j}\tilde{K}_{i},
\end{equation}
\begin{equation}
    K_{i}\tilde{\mathcal{E}}_{j} = (-1)^{a_{ij}}\tilde{\mathcal{E}}_{j}K_{i},
\end{equation}
\begin{equation}
    K_{i}\tilde{\mathcal{F}}_{j} = (-1)^{a_{ij}}\tilde{\mathcal{F}}_{j}K_{i},
\end{equation}
\begin{equation}
    \tilde{K}_{i}\mathcal{E}_{j} = (-1)^{a_{ij}}\mathcal{E}_{j}\tilde{K}_{i},
\end{equation}
\begin{equation}
    \tilde{K}_{i}\mathcal{F}_{j} = (-1)^{a_{ij}}\mathcal{F}_{j}\tilde{K}_{i}.
\end{equation}
\end{cor}
$\proof$
All these relations follow from (4.10), (4.11) if we take $\imath p$, $\imath s$, $\imath t$ to be equal to $1$ and $b^{-2}$ in different combinations.
$\Box$

\begin{cor}
Let $a_{ij} = 0$, $1\le i,j\le r$. Then
\begin{equation}
    [\mathcal{E}_{i},\mathcal{E}_{j}] = [\tilde{\mathcal{E}}_{i},\tilde{\mathcal{E}}_{j}] = [\tilde{\mathcal{E}}_{i},\mathcal{E}_{j}] = [\mathcal{E}_{i},\tilde{\mathcal{E}}_{j}]= 0,
\end{equation}
\begin{equation}
    [\mathcal{F}_{i},\mathcal{F}_{j}] = [\tilde{\mathcal{F}}_{i},\tilde{\mathcal{F}}_{j}] = [\tilde{\mathcal{F}}_{i},\mathcal{F}_{j}] = [\mathcal{F}_{i},\tilde{\mathcal{F}}_{j}]= 0.
\end{equation}
\end{cor}
$\proof$
These relations are consequences of (4.12), (4.13) if we take $\imath s$, $\imath t$ to be equal to $1$ and $b^{-2}$ in different combinations.
$\Box$

\begin{cor}
Let $i\ne j$, $1\le i,j\le r$. Then
\begin{equation}
    [\mathcal{E}_{i},\mathcal{F}_{j}] = [\tilde{\mathcal{E}}_{i},\tilde{\mathcal{F}}_{j}] = [\tilde{\mathcal{E}}_{i},\mathcal{F}_{j}] = [\mathcal{E}_{i},\tilde{\mathcal{F}}_{j}] = 0.
\end{equation}
\end{cor}
$\proof$
These relations follow from (4.16) if we take $\imath s$, $\imath t$ to be equal to $1$ and $b^{-2}$ in different combinations.
$\Box$

\begin{cor}
Let $a_{ij} = -1$ and let $N$, $M$ be non-negative integers. Then the following identities hold
\begin{equation}
    \mathcal{E}_{i}^{N}\mathcal{E}_{j}^{M} = q^{NM}\sum\limits_{n=0}^{min(N,M)}(-1)^{n}q^{-\frac{n^{2}}{2}+n}
    \frac{\prod\limits_{k=1}^{N}(1-q^{-2k})\prod\limits_{k=1}^{M}(1-q^{-2k})}
    {\prod\limits_{k=1}^{N-n}(1-q^{-2k})\prod\limits_{k=1}^{M-n}(1-q^{-2k})\prod\limits_{k=1}^{n}(1-q^{-2k})}
    \mathcal{E}_{j}^{M-n}\mathcal{E}_{ij}^{n}\mathcal{E}_{i}^{N-n},
\end{equation}
\begin{equation}
    \mathcal{F}_{i}^{N}\mathcal{F}_{j}^{M} = q^{NM}\sum\limits_{n=0}^{min(N,M)}(-1)^{n}q^{-\frac{n^{2}}{2}+n}
    \frac{\prod\limits_{k=1}^{N}(1-q^{-2k})\prod\limits_{k=1}^{M}(1-q^{-2k})}
    {\prod\limits_{k=1}^{N-n}(1-q^{-2k})\prod\limits_{k=1}^{M-n}(1-q^{-2k})\prod\limits_{k=1}^{n}(1-q^{-2k})}
    \mathcal{F}_{j}^{M-n}\mathcal{F}_{ij}^{n}\mathcal{F}_{i}^{N-n}.
\end{equation}
Substituting $N = M = 1$, one recovers the definition of non-simple roots \cite{Ip}, eq.(6.18)
\begin{equation}
    \mathcal{E}_{ij} = \frac{q^{\frac{1}{2}}\mathcal{E}_{j}\mathcal{E}_{i}-q^{-\frac{1}{2}}\mathcal{E}_{i}\mathcal{E}_{j}}{q-q^{-1}},
\end{equation}
\begin{equation}
    \mathcal{F}_{ij} = \frac{q^{\frac{1}{2}}\mathcal{F}_{j}\mathcal{F}_{i}-q^{-\frac{1}{2}}\mathcal{F}_{i}\mathcal{F}_{j}}{q-q^{-1}}.
\end{equation}
Substituting $N = 2$, $M = 1$ together with the definition of non-simple roots one recovers $q$- deformed Serre relations
\begin{equation}
    \mathcal{E}_{i}^{2}\mathcal{E}_{j} - (q+q^{-1})\mathcal{E}_{i}\mathcal{E}_{j}\mathcal{E}_{i} + \mathcal{E}_{j}\mathcal{E}_{i}^{2} = 0,
\end{equation}
for $a_{ij} = -1$, $1\le i,j\le r$. Analogous relation is true for $\mathcal{F}_{i}$, $\mathcal{F}_{j}$.
\end{cor}
$\proof$

$$
\mathcal{E}_{i}^{\imath s}\mathcal{E}_{j}^{\imath t} = e^{-\pi\imath b^{2}st}\int dz e^{\frac{\pi\imath b^{2}z^{2}}{2}-\pi bQz}
\frac{G_{b}(-\imath bs+\imath bz)G_{b}(-\imath bz)}{G_{b}(-\imath bs)}\frac{G_{b}(-\imath bt+\imath bz)}{G_{b}(-\imath bt)}
\mathcal{E}_{j}^{\imath t-\imath z}\mathcal{E}_{ij}^{\imath z}\mathcal{E}_{i}^{\imath s-\imath z}.
$$
Let $\imath s = N$. Using delta distribution formula
$$
    \frac{G_{b}(x)G_{b}(-Nb-x)}{G_{b}(-Nb)} =
    \sum\limits_{n=0}^{N}\frac{\prod\limits_{k=1}^{N}(1-q^{-2k})}{\prod\limits_{k=1}^{n}(1-q^{-2k})\prod\limits_{k=1}^{N-n}(1-q^{-2k})}
    \delta(x+nb),
$$
we evaluate the integral
$$
\mathcal{E}_{i}^{N}\mathcal{E}_{j}^{\imath t} = e^{-\pi b^{2}Nt}\int dz
e^{\frac{\pi\imath b^{2}z^{2}}{2}-\pi bQz}
\frac{G_{b}(-Nb+\imath bz)G_{b}(-\imath bz)}{G_{b}(-Nb)}\frac{G_{b}(-\imath bt+\imath bz)}{G_{b}(-\imath bt)}
\mathcal{E}_{j}^{\imath t-\imath z}\mathcal{E}_{ij}^{\imath z}\mathcal{E}_{i}^{N-\imath z} =
$$
$$
e^{-\pi b^{2}Nt}\sum\limits_{n=0}^{N}\frac{\prod\limits_{k=1}^{N}(1-q^{-2k})}{\prod\limits_{k=1}^{n}(1-q^{-2k})\prod\limits_{k=1}^{N-n}(1-q^{-2k})}
\int dz e^{\frac{\pi\imath b^{2}z^{2}}{2}-\pi bQz}\delta(-\imath bz+nb)\frac{G_{b}(-\imath bt+\imath bz)}{G_{b}(-\imath bt)}
\mathcal{E}_{j}^{\imath t-\imath z}\mathcal{E}_{ij}^{\imath z}\mathcal{E}_{i}^{N-\imath z} =
$$
$$
e^{-\pi b^{2}Nt}\sum\limits_{n=0}^{N}(-1)^{n}q^{-\frac{n^{2}}{2}+n}
\frac{\prod\limits_{k=1}^{N}(1-q^{-2k})}{\prod\limits_{k=1}^{n}(1-q^{-2k})\prod\limits_{k=1}^{N-n}(1-q^{-2k})}
\frac{G_{b}(-\imath bt+nb)}{G_{b}(-\imath bt)}\mathcal{E}_{j}^{\imath t-n}\mathcal{E}_{ij}^{n}\mathcal{E}_{i}^{N-n}.
$$
Using functional equation for $G_{b}(z)$
\begin{equation}
\frac{G_{b}(x+nb)}{G_{b}(x)} =
\prod\limits_{k =0}^{n-1}(1-q^{2k}e^{2\pi\imath bx}),
\end{equation}
we obtain

$$
\mathcal{E}_{i}^{N}\mathcal{E}_{j}^{\imath t} =
e^{-\pi b^{2}Nt}\sum\limits_{n=0}^{N}(-1)^{n}q^{-\frac{n^{2}}{2}+n}
\frac{\prod\limits_{k=1}^{N}(1-q^{-2k})\prod\limits_{k=0}^{n-1}(1-q^{2k}e^{2\pi b^{2}t})}{\prod\limits_{k=1}^{n}(1-q^{-2k})\prod\limits_{k=1}^{N-n}(1-q^{-2k})}
\mathcal{E}_{j}^{\imath t-n}\mathcal{E}_{ij}^{n}\mathcal{E}_{i}^{N-n}.
$$
By taking $\imath t = M$ we prove the corollary
$$
\mathcal{E}_{i}^{N}\mathcal{E}_{j}^{M} = q^{NM}
\sum\limits_{n=0}^{N}(-1)^{n}q^{-\frac{n^{2}}{2}+n}
\frac{\prod\limits_{k=1}^{N}(1-q^{-2k})\prod\limits_{k=0}^{n-1}(1-q^{-2(M-k)})}
{\prod\limits_{k=1}^{n}(1-q^{-2k})\prod\limits_{k=1}^{N-n}(1-q^{-2k})}
\mathcal{E}_{j}^{M-n}\mathcal{E}_{ij}^{n}\mathcal{E}_{i}^{N-n} =
$$
$$
q^{NM}\sum\limits_{n=0}^{min(N,M)}(-1)^{n}q^{-\frac{n^{2}}{2}+n}
    \frac{\prod\limits_{k=1}^{N}(1-q^{-2k})\prod\limits_{k=1}^{M}(1-q^{-2k})}
    {\prod\limits_{k=1}^{N-n}(1-q^{-2k})\prod\limits_{k=1}^{M-n}(1-q^{-2k})\prod\limits_{k=1}^{n}(1-q^{-2k})}
    \mathcal{E}_{j}^{M-n}\mathcal{E}_{ij}^{n}\mathcal{E}_{i}^{N-n}.
$$
Now, let $\imath t = b^{-2}$, $N = 1$ to prove the second corollary
$$
\mathcal{E}_{i}\mathcal{E}_{j}^{b^{-2}} = -(\mathcal{E}_{j}^{b^{-2}}\mathcal{E}_{i} - q^{-\frac{1}{2}+1}(1-e^{-2\pi\imath})\mathcal{E}_{j}^{b^{-2}-1}\mathcal{E}_{ij}) =
-\mathcal{E}_{j}^{b^{-2}}\mathcal{E}_{i}.
$$
Analogously one can prove the same relations for $\mathcal{F}_{j}$, $1\le j\le r$.

$\Box$

In the proof of the previous corollary we have also obtained the following
\begin{cor}
Let $a_{ij} = -1$ and let $\tilde{\mathcal{E}}_{j} = \mathcal{E}_{j}^{b^{-2}}$, $\tilde{\mathcal{F}}_{j} = \mathcal{F}_{j}^{b^{-2}}$, $1\le j\le r$. Then the following identities hold
\begin{equation}
    \mathcal{E}_{i}\tilde{\mathcal{E}}_{j} = -\tilde{\mathcal{E}}_{j}\mathcal{E}_{i},
\end{equation}
\begin{equation}
    \mathcal{F}_{i}\tilde{\mathcal{F}}_{j} = -\tilde{\mathcal{F}}_{j}\mathcal{F}_{i}.
\end{equation}
\end{cor}

\begin{cor}
Let $a_{ij} = -1$ and let $\tilde{\mathcal{E}}_{j} = \mathcal{E}_{j}^{b^{-2}}$, $\tilde{\mathcal{F}}_{j} = \mathcal{F}_{j}^{b^{-2}}$, $\tilde{\mathcal{E}}_{ij} = \mathcal{E}_{ij}^{b^{-2}}$, $\tilde{\mathcal{F}}_{ij} = \mathcal{F}_{ij}^{b^{-2}}$, $1\le i,j\le r$. Then the following identities hold
\begin{equation}
    \tilde{\mathcal{E}}_{i}^{N}\tilde{\mathcal{E}}_{j}^{M} = \tilde{q}^{NM}\sum\limits_{n=0}^{min(N,M)}(-1)^{n}\tilde{q}^{-\frac{n^{2}}{2}+n}
    \frac{\prod\limits_{k=1}^{N}(1-\tilde{q}^{-2k})\prod\limits_{k=1}^{M}(1-\tilde{q}^{-2k})}
    {\prod\limits_{k=1}^{N-n}(1-\tilde{q}^{-2k})\prod\limits_{k=1}^{M-n}(1-\tilde{q}^{-2k})\prod\limits_{k=1}^{n}(1-\tilde{q}^{-2k})}
    \tilde{\mathcal{E}}_{j}^{M-n}\tilde{\mathcal{E}}_{ij}^{n}\tilde{\mathcal{E}}_{i}^{N-n},
\end{equation}
\begin{equation}
    \tilde{\mathcal{F}}_{i}^{N}\tilde{\mathcal{F}}_{j}^{M} = \tilde{q}^{NM}\sum\limits_{n=0}^{min(N,M)}(-1)^{n}\tilde{q}^{-\frac{n^{2}}{2}+n}
    \frac{\prod\limits_{k=1}^{N}(1-\tilde{q}^{-2k})\prod\limits_{k=1}^{M}(1-\tilde{q}^{-2k})}
    {\prod\limits_{k=1}^{N-n}(1-\tilde{q}^{-2k})\prod\limits_{k=1}^{M-n}(1-\tilde{q}^{-2k})\prod\limits_{k=1}^{n}(1-\tilde{q}^{-2k})}
    \tilde{\mathcal{F}}_{j}^{M-n}\tilde{\mathcal{F}}_{ij}^{n}\tilde{\mathcal{F}}_{i}^{N-n}.
\end{equation}
Taking in these formulas $N = M = 1$ one obtains the following expressions for the dual non-simple roots
\begin{equation}
    \tilde{\mathcal{E}}_{ij} = \frac{\tilde{q}^{\frac{1}{2}}\tilde{\mathcal{E}}_{j}\tilde{\mathcal{E}}_{i}-\tilde{q}^{-\frac{1}{2}}\tilde{\mathcal{E}}_{i}\tilde{\mathcal{E}}_{j}}{\tilde{q}-\tilde{q}^{-1}},
\end{equation}
\begin{equation}
    \tilde{\mathcal{F}}_{ij} = \frac{\tilde{q}^{\frac{1}{2}}\tilde{\mathcal{F}}_{j}\tilde{\mathcal{F}}_{i}-\tilde{q}^{-\frac{1}{2}}\tilde{\mathcal{F}}_{i}\tilde{\mathcal{F}}_{j}}{\tilde{q}-\tilde{q}^{-1}}.
\end{equation}
Taking $N = 2$, $M = 1$ and using these expressions one obtains the dual $q$-deformed Serre relations
\begin{equation}
    \tilde{\mathcal{E}}_{i}^{2}\tilde{\mathcal{E}}_{j} - (\tilde{q}+\tilde{q}^{-1})\tilde{\mathcal{E}}_{i}\tilde{\mathcal{E}}_{j}\tilde{\mathcal{E}}_{i} + \tilde{\mathcal{E}}_{j}\tilde{\mathcal{E}}_{i}^{2} = 0,
\end{equation}
for $a_{ij} = -1$, $1\le i,j\le r$. Similar relation is valid for $\tilde{\mathcal{F}}_{i}$, $\tilde{\mathcal{F}}_{j}$.
\end{cor}
$\proof$
The proof is similar to the proof of analogous relation for untilded generators.
$\Box$

\begin{cor}
The following formulas for coproduct hold
\begin{equation}
\Delta K_{j} = K_{j}\otimes K_{j},
\end{equation}
\begin{equation}
\Delta \mathcal{E}_{j} = \mathcal{E}_{j}\otimes 1 + K_{j}^{-1}\otimes \mathcal{E}_{j},
\end{equation}
\begin{equation}
\Delta \mathcal{F}_{j} = 1\otimes \mathcal{F}_{j} + \mathcal{F}_{j}\otimes K_{j},
\end{equation}
\begin{equation}
\Delta \tilde{K}_{j} = \tilde{K}_{j}\otimes \tilde{K}_{j},
\end{equation}
\begin{equation}
\Delta \tilde{\mathcal{E}}_{j} = \tilde{\mathcal{E}}_{j}\otimes 1 + \tilde{K}_{j}^{-1}\otimes \tilde{\mathcal{E}}_{j},
\end{equation}
\begin{equation}
\Delta \tilde{\mathcal{F}}_{j} = 1\otimes \tilde{\mathcal{F}}_{j} + \tilde{\mathcal{F}}_{j}\otimes \tilde{K}_{j}.
\end{equation}
\end{cor}
$\proof$
Coproduct for generators $K_{j}$, $1\le j\le r$ are trivial consequences of the formula (4.18) for $\imath p$ equal to $1$ and $b^{-2}$. The rest expressions follow from (4.19), (4.20) for $\imath s$, $\imath t$ equal to $1$ and $b^{-2}$ by evaluating the integral using the delta distribution formula (\ref{delta}).
$\Box$

\section{Modular double $M_{q\tilde{q}}(\mathfrak{g})$}

Let $q = e^{\pi\imath b^{2}}$, $\tilde{q} = e^{\pi\imath b^{-2}}$, $b^{2}\in \mathbb{R}\setminus\mathbb{Q}$ and let $K_{j}$, $E_{j}$, $F_{j}$, $1\le j\le r$ be generators of quantum group $U_{q}(\mathfrak{g})$. The rescaled generators are defined by
\begin{equation}
    \mathcal{E}_{j} = -\imath (q-q^{-1})E_{j},
\end{equation}
\begin{equation}
    \mathcal{F}_{j} = -\imath (q-q^{-1})F_{j}.
\end{equation}
The second set of generators $\tilde{K}_{j}$, $\tilde{E}_{j}$, $\tilde{F}_{j}$, $1\le j\le r$ is defined as follows
\begin{equation}\label{first transcendental relation}
    \tilde{K}_{j} = K_{j}^{b^{-2}},
\end{equation}
\begin{equation}
    \tilde{E}_{j} = \frac{\imath}{\tilde{q}-\tilde{q}^{-1}}\mathcal{E}_{j}^{b^{-2}},
\end{equation}
\begin{equation}\label{last transcendental relation}
    \tilde{F}_{j} = \frac{\imath}{\tilde{q}-\tilde{q}^{-1}}\mathcal{F}_{j}^{b^{-2}}.
\end{equation}

\begin{te}
The elements $K_{i}$, $E_{j}$, $F_{j}$ and $\tilde{K}_{i}$, $\tilde{E}_{j}$, $\tilde{F}_{j}$, $1\le j\le r$ generate a Hopf subalgebra of $A(\mathfrak{g})$ with the first set of generators satisfying the relations of $U_{q}(\mathfrak{g})$ and the second set of generators satisfying the relations of $U_{\tilde{q}}(\mathfrak{g})$ and both sets of generators satisfying the cross-relations:
\begin{equation}\label{first cross-relation}
    K_{i}\tilde{K}_{j} = \tilde{K}_{j}K_{i},
\end{equation}
\begin{equation}
    K_{i}\tilde{E}_{j} = (-1)^{a_{ij}}\tilde{E}_{j}K_{i},
\end{equation}
\begin{equation}
    \tilde{K}_{i}E_{j} = (-1)^{a_{ij}}E_{j}\tilde{K}_{i},
\end{equation}
\begin{equation}
    K_{i}\tilde{F}_{j} = (-1)^{a_{ij}}\tilde{F}_{j}K_{i},
\end{equation}
\begin{equation}
    \tilde{K}_{i}F_{j} = (-1)^{a_{ij}}F_{j}\tilde{K}_{i},
\end{equation}
\begin{equation}
    E_{i}\tilde{E}_{j} = (-1)^{a_{ij}}\tilde{E}_{j}E_{i},
\end{equation}
\begin{equation}
    F_{i}\tilde{F}_{j} = (-1)^{a_{ij}}\tilde{F}_{j}F_{i},
\end{equation}
\begin{equation}
    E_{i}\tilde{F}_{j} = \tilde{F}_{j}E_{i},
\end{equation}
\begin{equation}\label{last cross-relation}
    \tilde{E}_{i}F_{j} = F_{j}\tilde{E}_{i},
\end{equation}
where $1\le i,j\le r$.
\end{te}
$\proof$
The results of the Corollaries 3.1-3.3 and Corollaries 4.1-4.9 give the statement of the theorem.
$\Box$

 The defining relations of modular double have been observed in representations in \cite{F2} for the case of $U_{q}(\mathfrak{sl}(2))$, in \cite{GKL2} for the case of $U_{q}(\mathfrak{gl}(N))$ and in \cite{Ip2} for the case of other Lie algebras $\mathfrak{g}$.

\begin{de}
Modular double $M_{q\tilde{q}}(\mathfrak{g})$ is a Hopf algebra generated by generators $K_{i}$, $E_{j}$, $F_{j}$, $1\le j\le r$ satisfying the relations of $U_{q}(\mathfrak{g})$ and generators $\tilde{K}_{i}$, $\tilde{E}_{j}$, $\tilde{F}_{j}$, $1\le j\le r$ satisfying the relations of $U_{\tilde{q}}(\mathfrak{g})$ with the cross-relations (\ref{first cross-relation}-\ref{last cross-relation}). The transcendental relations (\ref{first transcendental relation}-\ref{last transcendental relation}) are not imposed.
\end{de}
 Note, that in this algebraic definition of modular double we do not require transcendental relations between two sets of generators. The modular double defined in such a way has two types of representations. Below we give an example of the representation of the modular double $M_{q\tilde{q}}(\mathfrak{gl}(N))$ for which there are no transcendental relations. \\*
Let us introduce more general parametrization $q = e^{\frac{\pi\imath\omega_{1}}{\omega_{2}}}$, $\tilde{q} = e^{\frac{\pi\imath\omega_{2}}{\omega_{1}}}$. This notation is reduced to the one we used in this paper if we set $\omega_{1} = b$, $\omega_{2} = b^{-1}$.
$U_{q}(\mathfrak{gl}(N))$ is generated by the elements $K_{n}$, $n = 1,...,N$ and $E_{n,n+1}$, $E_{n+1,n}$, $n = 1,...,N-1$ subjected to the following set of relations
\begin{equation}
E_{n,n+1}E_{m+1,m} -E_{m+1,m}E_{n,n+1} =
\delta_{nm}\frac{K_{n}K_{n+1}^{-1}-K_{n}^{-1}K_{n+1}}{q-q^{-1}},
\end{equation}
\begin{equation}
K_{n}E_{m,m+1} = q^{\delta_{nm}-\delta_{n,m+1}}E_{m,m+1}K_{n},
\end{equation}
\begin{equation}
K_{n}E_{m+1,m} =
q^{\delta_{n,m+1}-\delta_{nm}}E_{m+1,m}K_{n},
\end{equation}
together with quantum analogues of Serre relations
\begin{equation}
E_{n,n+1}E_{m,m+1} - E_{m,m+1}E_{n,n+1} =0, m\ne n\pm 1,
\end{equation}
\begin{equation}
E_{n,n+1}^{2}E_{n+1,n+2} -
(q+q^{-1})E_{n,n+1}E_{n+1,n+2}E_{n,n+1}+E_{n+1,n+2}E_{n,n+1}^{2} =0,
\end{equation}
\begin{equation}
E_{n+1,n+2}^{2}E_{n,n+1} -
(q+q^{-1})E_{n+1,n+2}E_{n,n+1}E_{n+1,n+2}+E_{n,n+1}E_{n+1,n+2}^{2} =0,
\end{equation}
\begin{equation}
E_{n+1,n}E_{m+1,m} - E_{m+1,m}E_{n+1,n} =0, m\ne n\pm 1,
\end{equation}
\begin{equation}
E_{n+1,n}^{2}E_{n+2,n+1} -
(q+q^{-1})E_{n+1,n}E_{n+2,n+1}E_{n+1,n}+E_{n+2,n+1}E_{n+1,n}^{2} =0,
\end{equation}
\begin{equation}
E_{n+2,n+1}^{2}E_{n+1,n} -
(q+q^{-1})E_{n+2,n+1}E_{n+1,n}E_{n+2,n+1}+E_{n+1,n}E_{n+2,n+1}^{2} =0.
\end{equation}
Analogously, the dual quantum group $U_{\tilde{q}}(\mathfrak{gl}(N))$ is generated by the elements $\tilde{K}_{n}$, $n = 1,...,N$ and $\tilde{E}_{n,n+1}$, $\tilde{E}_{n+1,n}$, $n = 1,...,N-1$ subjected to the relations
\begin{equation}
\tilde{E}_{n,n+1}\tilde{E}_{m+1,m} -\tilde{E}_{m+1,m}\tilde{E}_{n,n+1} =
\delta_{nm}\frac{\tilde{K}_{n}\tilde{K}_{n+1}^{-1}-\tilde{K}_{n}^{-1}\tilde{K}_{n+1}}{\tilde{q}-\tilde{q}^{-1}},
\end{equation}
\begin{equation}
\tilde{K}_{n}\tilde{E}_{m,m+1} = \tilde{q}^{\delta_{nm}-\delta_{n,m+1}}\tilde{E}_{m,m+1}\tilde{K}_{n},
\end{equation}
\begin{equation}
\tilde{K}_{n}\tilde{E}_{m+1,m} =
\tilde{q}^{\delta_{n,m+1}-\delta_{nm}}\tilde{E}_{m+1,m}\tilde{K}_{n},
\end{equation}
together with quantum analogues of Serre relations
\begin{equation}
\tilde{E}_{n,n+1}\tilde{E}_{m,m+1} - \tilde{E}_{m,m+1}\tilde{E}_{n,n+1} =0, m\ne n\pm 1,
\end{equation}
\begin{equation}
\tilde{E}_{n,n+1}^{2}\tilde{E}_{n+1,n+2} -
(\tilde{q}+\tilde{q}^{-1})\tilde{E}_{n,n+1}\tilde{E}_{n+1,n+2}\tilde{E}_{n,n+1}+\tilde{E}_{n+1,n+2}\tilde{E}_{n,n+1}^{2} =0,
\end{equation}
\begin{equation}
\tilde{E}_{n+1,n+2}^{2}\tilde{E}_{n,n+1} -
(\tilde{q}+\tilde{q}^{-1})\tilde{E}_{n+1,n+2}\tilde{E}_{n,n+1}\tilde{E}_{n+1,n+2}+\tilde{E}_{n,n+1}\tilde{E}_{n+1,n+2}^{2} =0,
\end{equation}
\begin{equation}
\tilde{E}_{n+1,n}\tilde{E}_{m+1,m} - \tilde{E}_{m+1,m}\tilde{E}_{n+1,n} =0, m\ne n\pm 1,
\end{equation}
\begin{equation}
\tilde{E}_{n+1,n}^{2}\tilde{E}_{n+2,n+1} -
(\tilde{q}+\tilde{q}^{-1})\tilde{E}_{n+1,n}\tilde{E}_{n+2,n+1}\tilde{E}_{n+1,n}+\tilde{E}_{n+2,n+1}\tilde{E}_{n+1,n}^{2} =0,
\end{equation}
\begin{equation}
\tilde{E}_{n+2,n+1}^{2}\tilde{E}_{n+1,n} -
(\tilde{q}+\tilde{q}^{-1})\tilde{E}_{n+2,n+1}\tilde{E}_{n+1,n}\tilde{E}_{n+2,n+1}+\tilde{E}_{n+1,n}\tilde{E}_{n+2,n+1}^{2} =0.
\end{equation}
The modular double $M_{q\tilde{q}}(\mathfrak{gl}(N))$ is generated by the generators of $U_{q}(\mathfrak{gl}(N))$ and generators of $U_{\tilde{q}}(\mathfrak{gl}(N))$ with the following cross-relations:
\begin{equation}
E_{n,n+1}\tilde{K}_{m} = (-1)^{\delta_{nm}+\delta_{n,m-1}}\tilde{K}_{m}E_{n,n+1},
\end{equation}
\begin{equation}
E_{n+1,n}\tilde{K}_{m} = (-1)^{\delta_{nm}+\delta_{n,m-1}}\tilde{K}_{m}E_{n+1,n},
\end{equation}
\begin{equation}
E_{n,n+1}\tilde{E}_{m,m+1} = (-1)^{\delta_{n,m+1}+\delta_{n+1,m}}\tilde{E}_{m,m+1}E_{n,n+1},
\end{equation}
\begin{equation}
E_{n,n+1}\tilde{E}_{m+1,m} = \tilde{E}_{m+1,m}E_{n,n+1},
\end{equation}
\begin{equation}
E_{n+1,n}\tilde{E}_{m+1,m} = (-1)^{\delta_{n,m-1}+\delta_{n-1,m}}\tilde{E}_{m+1,m}E_{n+1,n},
\end{equation}
\begin{equation}
\tilde{E}_{n,n+1}K_{m} = (-1)^{\delta_{nm}+\delta_{n,m-1}}K_{m}\tilde{E}_{n,n+1},
\end{equation}
\begin{equation}
\tilde{E}_{n+1,n}K_{m} = (-1)^{\delta_{nm}+\delta_{n,m-1}}K_{m}\tilde{E}_{n+1,n},
\end{equation}
\begin{equation}
\tilde{E}_{n,n+1}E_{m+1,m} = E_{m+1,m}\tilde{E}_{n,n+1}.
\end{equation}

\begin{te}\cite{GKL2}
Let $q = e^{\frac{\pi\imath\omega_{1}}{\omega_{2}}}$ and $\tilde{q} = e^{\frac{\pi\imath\omega_{2}}{\omega_{1}}}$. The following operators define a representation   of modular double $M_{q\tilde{q}}(\mathfrak{gl}(N))$
\begin{equation}
E_{n,n+1} = \frac{2e^{\frac{\pi\imath(\omega_{1}+\omega_{2})(n-1)}{2\omega_{2}}}}{q-q^{-1}}
\sum\limits_{j = 1}^{n}
\frac{\prod\limits_{r = 1}^{n+1}\sinh\frac{\pi}{\omega_{2}}(\gamma_{nj}-\gamma_{n+1,r}-\frac{\imath}{2}(\omega_{1}+\omega_{2}))}
{\prod\limits_{s \ne j}^{n}\sinh\frac{\pi}{\omega_{2}}(\gamma_{nj}-\gamma_{ns})}
e^{-\imath\omega_{1}\partial_{\gamma_{nj}}},
\end{equation}
\begin{equation}
E_{n+1,n} =
\frac{2e^{-\frac{\pi\imath(\omega_{1}+\omega_{2})(n-1)}{2\omega_{2}}}}{q-q^{-1}}
\sum\limits_{j=1}^{n}
\frac{\prod\limits_{r=1}^{n-1}
\sinh\frac{\pi}{\omega_{2}}(\gamma_{nj}-\gamma_{n-1,r}+\frac{\imath}{2}(\omega_{1}+\omega_{2}))}
{\prod\limits_{s\ne j}^{n}\sinh\frac{\pi}{\omega_{2}}(\gamma_{nj}-\gamma_{ns})}
e^{\imath\omega_{1}\partial_{\gamma_{nj}}},
\end{equation}
\begin{equation}
K_{n} = e^{\frac{\pi}{\omega_{2}}(\sum\limits_{j=1}^{n}\gamma_{nj}-\sum\limits_{j=1}^{n-1}\gamma_{n-1,j})},
\end{equation}

\begin{equation}
\tilde{E}_{n,n+1} = \frac{2e^{\frac{\pi\imath(\omega_{1}+\omega_{2})(n-1)}{2\omega_{1}}}}{\tilde{q}-\tilde{q}^{-1}}
\sum\limits_{j = 1}^{n}
\frac{\prod\limits_{r = 1}^{n+1}\sinh\frac{\pi}{\omega_{1}}(\gamma_{nj}-\gamma_{n+1,r}-\frac{\imath}{2}(\omega_{1}+\omega_{2}))}
{\prod\limits_{s \ne j}^{n}\sinh\frac{\pi}{\omega_{1}}(\gamma_{nj}-\gamma_{ns})}
e^{-\imath\omega_{2}\partial_{\gamma_{nj}}},
\end{equation}
\begin{equation}
\tilde{E}_{n+1,n} =
\frac{2e^{-\frac{\pi\imath(\omega_{1}+\omega_{2})(n-1)}{2\omega_{1}}}}{\tilde{q}-\tilde{q}^{-1}}
\sum\limits_{j=1}^{n}
\frac{\prod\limits_{r=1}^{n-1}
\sinh\frac{\pi}{\omega_{1}}(\gamma_{nj}-\gamma_{n-1,r}+\frac{\imath}{2}(\omega_{1}+\omega_{2}))}
{\prod\limits_{s\ne j}^{n}\sinh\frac{\pi}{\omega_{1}}(\gamma_{nj}-\gamma_{ns})}
e^{\imath\omega_{2}\partial_{\gamma_{nj}}},
\end{equation}
\begin{equation}
\tilde{K}_{n} = e^{\frac{\pi}{\omega_{1}}(\sum\limits_{j=1}^{n}\gamma_{nj}-\sum\limits_{j=1}^{n-1}\gamma_{n-1,j})}.
\end{equation}
\end{te}
$\proof$ This is the same representation as was introduced in \cite{GKL2}, (sections 3.2-3.3) if one replaces $2\pi$ by $\pi$.
After this replacement it is simple to check that the cross-relations instead of commutativity relations  appear. $\Box$


There is another type of realization of the principal series representations \cite{FrIp}, \cite{Ip2} which are $q$-analogues of principal series representations of universal enveloping algebra of semisimple Lie algebra $\mathfrak{g}$ in Lusztig's parametrization. The classical limit of these representations was introduced earlier in \cite{GLO} in sections 2.4.1-2.4.4 for classical series of Lie algebras.

\section{Appendix}

\subsection{Quantum dilogarithm and its properties}
The basic properties of non-compact quantum dilogarithm/double sine listed below are extracted mainly from \cite{KLSTS}, \cite{BT}, \cite{V}.
Introduce the following notation $q = e^{\pi\imath b^{2}}$, $\tilde{q} = e^{\pi\imath b^{-2}}$, $Q = b+b^{-1}$, $\zeta_{b} = e^{\frac{\pi\imath}{4} + \frac{\pi\imath(b^{2}+b^{-2})}{12}}$.\\*
\\*
\textbf{The integral representation of $G_{b}(z)$:}
\begin{equation}
\log G_{b}(z) = \log\bar{\zeta}_{b} - \int\limits_{\mathbb{R}+\imath 0} \frac{dt}{t}\frac{e^{zt}}{(1-e^{bt})(1-e^{b^{-1}t})}.
\end{equation}
\textbf{Noncompact analog of $q$-exponential $g_{b}(z)$:}
\begin{equation}
g_{b}(z) = \frac{\bar{\zeta}_{b}}{G_{b}(\frac{Q}{2}+\frac{1}{2\pi\imath b}\log z)}.
\end{equation}
\textbf{Product representation:}
\begin{equation}
G_{b}(x) = \bar{\zeta}_{b}\frac{\prod\limits_{n=1}^{\infty}(1-e^{2\pi\imath b^{-1}(x-nb^{-1})})}{\prod\limits_{n=0}^{\infty}(1-e^{2\pi\imath b(x+nb)})},
\end{equation}
\begin{equation}
g_{b}(x) = \frac{\prod\limits_{n=0}^{\infty}(1+xq^{2n+1})}{\prod\limits_{n=0}^{\infty}(1+x^{b^{-2}}\tilde{q}^{-2n-1})}.
\end{equation}
\textbf{Functional equations:}
\begin{equation}
G_{b}(x +b^{\pm 1}) = (1-e^{2\pi\imath b^{\pm 1}x})G_{b}(x),
\end{equation}
or more generally
\begin{equation} \label{func eq}
\frac{G_{b}(x+n_{1}b+n_{2}b^{-1})}{G_{b}(x)} =
\prod\limits_{k_{1} =0}^{n_{1}-1}(1-q^{2k_{1}}e^{2\pi\imath bx})\prod\limits_{k_{2} =0}^{n_{2}-1}(1-\tilde{q}^{2k_{2}}e^{2\pi\imath b^{-1}x}),
\end{equation}
\begin{equation}
g_{b}(q^{-1}x) = (1+x)g_{b}(qx).
\end{equation}
\textbf{Reflection formula:}
\begin{equation}\label{reflection}
G_{b}(x)G_{b}(Q-x) = e^{\pi\imath x(x-Q)}.
\end{equation}
\textbf{Poles and zeros:}
\begin{equation}\begin{split}\label{Poles}
\lim_{x\rightarrow 0}xG_{b}(x-n_{1}b-n_{2}b^{-1}) = \frac{1}{2\pi}\prod\limits_{k_{1}=1}^{n_{1}}(1-q^{-2k_{1}})^{-1}\prod\limits_{k_{2}=1}^{n_{2}}(1-\tilde{q}^{-2k_{2}})^{-1},    \\
\lim_{x\rightarrow 0}xG_{b}^{-1}(x+Q+n_{1}b+n_{2}b^{-1}) =      \\
\frac{1}{2\pi}(-1)^{n_{1}+n_{2}+1}q^{-n_{1}(n_{1}+1)}\tilde{q}^{-n_{2}(n_{2}+1)}\prod\limits_{k_{1}=1}^{n_{1}}(1-q^{-2k_{1}})^{-1}\prod\limits_{k_{2}=1}^{n_{2}}(1-\tilde{q}^{-2k_{2}})^{-1}.
\end{split}\end{equation}
\textbf{Tau-binomial integral} \cite{FKV},\cite{Ka},\cite{PT}:
\begin{equation}\label{tau-integral}
\int\limits_{\mathcal{C}} d\tau e^{-2\pi b\beta\tau}\frac{G_{b}(\alpha+\imath b\tau)}{G_{b}(Q+\imath b\tau)} =
\frac{G_{b}(\alpha)G_{b}(\beta)}{G_{b}(\alpha+\beta)},
\end{equation}
where the contour $\mathcal{C}$ goes along the real axis above the sequences of poles going down and below sequences of poles going up.\\*
\\*
\textbf{Delta distributions \cite{Ip}:}
\begin{equation}\begin{split}  \label{delta}
    \frac{G_{b}(x)G_{b}(-N_{1}b-N_{2}b^{-1}-x)}{G_{b}(-N_{1}b-N_{2}b^{-1})} =
    \sum\limits_{n_{1}=0}^{N_{1}}\sum\limits_{n_{2}=0}^{N_{2}}\frac{\prod\limits_{k_{1}=1}^{N_{1}}(1-q^{-2k_{1}})}{\prod\limits_{k_{1}=1}^{n_{1}}(1-q^{-2k_{1}})\prod\limits_{k_{1}=1}^{N_{1}-n_{1}}(1-q^{-2k_{1}})} \\
   \times \frac{\prod\limits_{k_{2}=1}^{N_{2}}(1-\tilde{q}^{-2k_{2}})}{\prod\limits_{k_{2}=1}^{n_{2}}(1-\tilde{q}^{-2k_{2}})\prod\limits_{k_{2}=1}^{N_{2}-n_{2}}(1-\tilde{q}^{-2k_{2}})}
    \delta(x+n_{1}b+n_{2}b^{-1}).
\end{split}\end{equation}
\\*
   \textbf{ 4-5 relation \cite{V}:}
    \begin{equation}\label{4-5}
        \int\limits_{\mathcal{C}} d\tau e^{2\pi\imath(\imath\alpha+\tau)(\imath\beta+\tau)}G_{b}(\alpha-\imath\tau)G_{b}(\beta-\imath\tau)G_{b}(\gamma+\imath\tau)G_{b}(\imath\tau) =
        \frac{G_{b}(\alpha)G_{b}(\beta)G_{b}(\alpha+\gamma)G_{b}(\beta+\gamma)}{G_{b}(\alpha+\beta+\gamma)},
    \end{equation}
    where the contour $\mathcal{C}$ goes along the real axis above the sequences of poles going down and below sequences of poles going up.\\*
    \\*
\textbf{6-9 identity \cite{V}:}
\begin{equation}\begin{split}        \label{6-9}
\frac{G_{b}(A)G_{b}(B)G_{b}(C)G_{b}(A+D)G_{b}(B+D)G_{b}(C+D)}
{G_{b}(A+B+D)G_{b}(A+C+D)G_{b}(B+C+D)} = \\
\int\limits_{\mathcal{C}} d\tau e^{2\pi\imath\tau^{2}-2\pi D\tau}
\frac{G_{b}(A+\imath\tau)G_{b}(B+\imath\tau)G_{b}(C+\imath\tau)G_{b}(D-\imath\tau)G_{b}(-\imath\tau)}{G_{b}(A+B+C+D+\imath\tau)},
\end{split}\end{equation}
where the contour $\mathcal{C}$ goes along the real axis above the sequences of poles going down and below sequences of poles going up.\\*
\\*
\textbf{$q$-binomial theorem \cite{BT}:}\\*
Let $u$, $v$ be positive self-adjoint operators subject to the relations $uv = q^{2}vu$.
Then:
\begin{equation}\label{q-binomial}
(u+v)^{\imath s} = \int\limits_{\mathcal{C}} d\tau \frac{G_{b}(-\imath b\tau)G_{b}(-\imath bs+\imath b\tau)}{G_{b}(-\imath bs)}u^{\imath s-\imath\tau}v^{\imath\tau},
\end{equation}
where the contour $\mathcal{C}$ goes along the real axis above the sequences of poles going down and below sequences of poles going up.

\end{document}